\documentclass{jT}

\usepackage{amssymb,amsmath,latexsym,amsthm}

\newtheorem{theorem}{Theorem}[section]

\theoremstyle{definition}
\newtheorem*{definition}{Definition}
\newtheorem*{remark}{Remark}

\newtheorem*{question}{Question}
\numberwithin{equation}{section}

\begin{document}

\title[Thue and the Thue-Morse sequence]{Thue, Combinatorics on words, and conjectures inspired by the Thue-Morse sequence}

\author{\sc Jean-Paul ALLOUCHE}
\address{Jean-Paul ALLOUCHE\\
Institut de Math\'ematiques de Jussieu\\
\'Equipe Combinatoire et Optimisation\\
Universit\'e Pierre et Marie Curie, Case 247\\
4 Place Jussieu\\
F-75252 Paris Cedex 05 France}
\email{allouche@math.jussieu.fr}
\urladdr{http://www.math.jussieu.fr/$\sim$allouche}

\subjclass[2010]{11 B 85, 68 R 15}

\maketitle

\begin{center}
On the occasion of Axel Thue's 150th birthday
\end{center}

\begin{resume}
Nous d\'ecrivons quelques r\'esultats r\'ecents sur la suite de Thue-Morse, ainsi que des
questions ou conjectures, dont l'une, due \`a Shevelev, est r\'esolue dans cet article.
\end{resume}

\begin{abstr}
We describe some recent results on the Thue-Morse sequence. We also list open questions and
conjectures, one of which is due to Shevelev and proved in this paper.
\end{abstr}

\bigskip

\section{Introduction}

The birth of combinatorics on words, i.e., the study of finite sequences (called {\em words}) taking 
their values in a finite set (called {\em alphabet}), aka the study of the free monoid generated 
by a finite set, can be dated to 1906, when the first of the two seminal papers \cite{Thue1, Thue2} 
of Axel Thue appeared. Thue was interested in constructing an infinite sequence on three symbols with no
squares (i.e., without two consecutive identical blocks) in it, and/or an infinite sequence on two symbols
with no cubes (i.e., without three consecutive identical blocks) in it. Thue's sequence on two symbols
was actually already described in a number-theoretic paper by Prouhet \cite{Prouhet} addressing what
is now known as the Prouhet-Tarry-Escott problem (see, e.g., \cite{BorIng}), so that to some extent 
Prouhet might be considered as one of the fathers of combinatorics of words (see \cite{Cerny}). 

\bigskip

To begin with we would like to point out the first few lines of Thue's second paper \cite{Thue2}, 
where he writes (translation by Berstel in \cite{Berstel}):

\textit{For the development of logical sciences it will be important, without consideration
for possible applications, to find large domains for speculation about difficult problems. 
In this paper, we present some investigations in the theory of sequences of symbols, a theory
that has some connections with number theory.}

\bigskip

Combinatorics on words developed extensively in the last forty years, interacting with many 
fields, from number theory to harmonic analysis, from theoretical computer science to physics,
with incursions towards tiling, music, etc.

\bigskip

We will briefly indicate some of these developments as well as pointers to surveys or to the original
bibliography. Then we will single out a few conjectures or open problems to show that the domain
is still full of really nontrivial and interesting questions.

\section{The Thue-Morse sequence and its numerous occurrences}

\subsection{Recalling the definition of the Thue-Morse sequence}

\begin{definition} The Thue-Morse sequence is the sequence ${\mathfrak t} = ({\mathfrak t}(n))_{n \geq 0}$
defined by ${\mathfrak t}(n) = 0$ if the sum of the binary digits of the integer $n$ is even, and
${\mathfrak t}(n) = 1$ if the sum of the binary digits of $n$ is odd. 
We let ${\mathfrak u} = ({\mathfrak u}(n))_{n \geq 0}$ denote the sequence defined by 
${\mathfrak u}(n) = (-1)^{{\mathfrak t}(n)}$.
\end{definition}

\begin{remark}
Several equivalent definitions can be given for the sequences ${\mathfrak t}$ and ${\mathfrak u}$.
In particular ${\mathfrak t}$ is the infinite fixed point, beginning with $0$, of the (uniform) morphism
$0 \to 01$, $1 \to 10$. For more about the sequences ${\mathfrak t}$ and ${\mathfrak u}$ we refer the 
reader to, e.g., \cite{ubiq}.
\end{remark}

\bigskip

After the papers of Thue \cite{Thue1, Thue2} on what is now called the Thue-Morse sequence (or the
Prouhet-Tue-Morse sequence), the sequence was rediscovered (in particular by Morse \cite{Morse})
and/or studied in very many different contexts. J. Shallit and the author wrote a survey some years ago 
entitled ``The ubiquitous Thue-Morse sequence'' \cite{ubiq}, but since that time, we have gathered more
and more new references on occurrences of this sequence in the literature. We would like to briefly 
allude here to only two ``new'' directions.

\subsection{Nonrepetitive coloring of graphs}

The seminal paper \cite{AGHR} introduced the notion of coloring graphs ``without repetitions'': 
a coloring of the set of edges of a graph $G$ is called non-repetitive if the sequence of colors 
on any path in $G$ is non-repetitive, i.e., has no squares (consecutive identical blocks) in it.
In the last ten years a very large number of papers have been devoted to this subject, and it would 
deserve a full survey by itself.

\subsection{The Thue-Morse sequence in games, economics, etc.}

Suppose that two players $A$ and $B$ are playing a sequential tournament where the first mover may 
have some advantage (or disadvantage) just by playing first. Then, if the players are to play
several consecutive games, it seems fair that the order be, e.g., $A \ B \ A \ B \ A \ B \ \ldots$
In order to try to improve the fairness of the tournament, one can imagine that the first player
plays, then each player --begining with the second one-- plays twice (this is for example the case 
for tiebreak in tennis), yielding the sequence $A \ B \ B \ A \ A \ B \ B \ A \ A \ B \ B \ A \ \ldots$
Trying to improve again and again on this alternation leads to... the Thue-Morse sequence: to cite but a 
few papers on this subject, see \cite{BraTay} cited in \cite{LevSta}; also see \cite{CooDut, KPR,Pal, Richman}.

\section{Some conjectures and open questions}

\subsection{A recent conjecture of Shevelev}

V. Shevelev wrote recently several papers on the Thue-Morse sequence or related problems, see 
\cite{Shev2007a}--\cite{Shevelev}, see also the two papers by Moses and Shevelev \cite{ShevMos1, ShevMos2}. 
These papers contain nice results, but also questions. One of Shevelev's questions about
the Thue-Morse sequence is rather intriguing \cite{Shevelev}. In order to give the flavor of this question, 
let us recall that an integer is called {\it evil} [in French {\it pa\"{\i}en}] (resp.\ {\it odious} 
[in French {\it impie}]) if the sum of its binary digits is even (resp.\ odd)
I learned in the paper of Fraenkel \cite{Fraenkel2} that the terminology, inspired phonetically
by the words ``even'' and ``odd'', was coined by the 
authors of \cite{BCK} when they were composing their book\footnote{In the papers \cite{Fraenkel1, Fraenkel2} 
the terms ``vile'' and ``dopey'' are proposed for other kinds of integers. We also cite a suggestion
of Ingrid Daubechies to call {\it perfidy} the fact of being evil or odious, as {\it parity} is 
the fact of being even or odd, see {\tt http://blog.tanyakhovanova.com/?p=97}}, see \cite[p.\ 463]{BCK2}. 
A particular case of Shevelev's question is: for any odd integer $a$, evil and odious numbers alternate 
in the increasing sequence of integers $n$ such that ${\mathfrak u}(n+a) = {\mathfrak u}(n)$ where 
${\mathfrak u}$ is the $\pm 1$ Thue-Morse sequence).

\medskip

More generally Shevelev made the following conjecture that we will prove here.

\begin{theorem}[Shevelev's conjecture]\label{shev}
Let $({\mathfrak u}(n))_{n \geq 0} = ((-1)^{{\mathfrak t}(n)})_{n \geq 0}$. Let $a$ be a positive integer.
Let $B_a := \{\ell_0 < \ell_1 < \ell_2 < \ldots\}$ (resp.\ $C_a = \{m_0 < m_1 < m_2 < \ldots\}$) 
denote the sequence of integers in increasing order satisfying ${\mathfrak u}(\ell + a) 
= - {\mathfrak u}(\ell)$ (resp.\ ${\mathfrak u}(m+a) = {\mathfrak u}(m)$).

\noindent
Let $(\beta_a(n))_{n \geq 0}$ (resp.\ $(\gamma_a(n))_{n \geq 0}$) be the $\pm 1$ 
sequences defined by
$$
\beta_a(n) := {\mathfrak u}(\ell_n), \ \gamma_a(n) := {\mathfrak u}(m_n).
$$
Then the sequences $\beta$ and $\gamma$ are periodic, of smallest period $2^{v(a)+1}$ where
$v(a)$ is the $2$-adic valuation of $a$ (i.e., the largest integer $k$ such that $2^k$
divides $a$). They satisfy $\beta = -\gamma$. Furthermore the prefix of length $2^{v(a)+1}$ 
of the sequence $\gamma_a$ is equal to the prefix of length $2^{v(a)+1}$ of the Thue-Morse 
sequence $u$ if ${\mathfrak u}(a) = 1$ and to ``minus this prefix'' (i.e., where all $+1$ are replaced 
by $-1$ and conversely) if ${\mathfrak u}(a)=-1$. 
\end{theorem}

\proof We first prove the periodicity assertion for sequence $\gamma$. Note that, in order
to prove the assertions ``equal to the prefix...'' and ``equal to minus the prefix...'', it 
suffices to prove the assertion ``equal to $\pm 1$ the prefix...'' thanks to the remark that the 
sequence $\gamma_a$ begins with $1$ when ${\mathfrak u}(a)=1$ and in $-1$ when ${\mathfrak u}(a)=-1$. 
Namely $({\mathfrak u}(0),{\mathfrak u}(1),{\mathfrak u}(2)) = (+1, -1, -1)$; also, since $a$ is odd, 
$a+1$ even, thus ${\mathfrak u}(a+1)$ and ${\mathfrak u}(a+2)$ have opposite signs: namely, if $a = 2b+1$,
then ${\mathfrak u}(a+1) = {\mathfrak u}(2b+2) = {\mathfrak u}(b+1)$, while ${\mathfrak u}(a+2) =
{\mathfrak u}(2b+3) = {\mathfrak u}(2(b+1)+1) = - {\mathfrak u}(b+1)$. Hence 
$({\mathfrak u}(a),{\mathfrak u}(a+1),{\mathfrak u}(a+2)) = ({\mathfrak u}(a), +1, -1)$ or 
$({\mathfrak u}(a), -1, +1)$. So, if ${\mathfrak u}(a)=+1$, we have $\min C_a = 0$, thus 
${\mathfrak u}(\min C_a) = {\mathfrak u}(0) = +1$. 
If ${\mathfrak u}(a)=-1$, we have $\min C_a = 2$ or $1$, thus ${\mathfrak u}(\min C_a) = -1$.

\medskip
\noindent
Then the proof splits into two parts depending on the parity of $a$.

\medskip

-- If $a$ is odd, then consider the sequence $(w_a(n))_{n \geq 0}$ where 
$w_a(n) := {\mathfrak u}(a+n) + {\mathfrak u}(n)$. Clearly $w_a(n) = \pm 2$ if $n$ belongs to $C_a$, and
$w_a(n) = 0$ otherwise. More precisely $w_a(n) = 2$ if $n$ belongs to $C_a$ and ${\mathfrak u}(n) = 1$, and 
$w_a(n) = -2$ if $n$ belongs to $C_a$ and ${\mathfrak u}(n) = -1$. In other words $(w_a(n))_{n \geq 0}$ is 
the sequence obtained by interspersing blocks of zeros into the sequence $(2\gamma_a(n))_{n \geq 0}$.
What we have to prove is that the sequence $(2\gamma_a(n))_{n \geq 0}$ is one of the alternating sequences
$(+2 \ -2)^{\infty}$ or $(-2 \ +2)^{\infty}$. Since $2\gamma_a$ takes only the values $+2$ and $-2$, this 
is equivalent to proving that the summatory function of sequence $w_a$ either takes only the values 
$0$, $+2$, or only the values $0$, $-2$. Now the summatory function of $w_a$ is the sequence of coefficients 
of the formal power series $S_a(X) := \frac{1}{1-X}\sum_{n \geq 0} w_a(n)X^n$. We have
$$
\begin{array}{lll}
S_a(X) &=& \displaystyle\frac{1}{1-X}\sum_{n \geq 0} ({\mathfrak u}(n+a) + {\mathfrak u}(n)) X^n \\
&=& \displaystyle\frac{1}{1-X} 
      \left(\sum_{m \geq a} {\mathfrak u}(m)X^{m-a} + \sum_{n \geq 0} {\mathfrak u}(n) X^n\right). 
\end{array}
$$
Hence, multiplying by $X^a$,
$$
X^aS_a(X) = 
- \frac{1}{1-X}\left(\sum_{0 \leq m \leq a-1} {\mathfrak u}(m) X^m\right) 
+ \frac{1+X^a}{1-X} \sum_{n \geq 0} {\mathfrak u}(n) X^n. \\
$$
But
$$
\begin{array}{lll}
\displaystyle\sum_{n \geq 0} {\mathfrak u}(n) X^n &=& 
\displaystyle\sum_{n \geq 0} {\mathfrak u}(2n) X^{2n} 
+ \displaystyle\sum_{n \geq 0} {\mathfrak u}(2n+1) X^{2n+1} \\
&=& \displaystyle\sum_{n \geq 0} {\mathfrak u}(n) X^{2n} 
    - \displaystyle\sum_{n \geq 0} {\mathfrak u}(n) X^{2n+1} \\
&=& \displaystyle(1-X)\sum_{n \geq 0} {\mathfrak u}(n) X^{2n}
\end{array}
$$
(this is well-known, as is the iteration yielding that $\sum_{n \geq 0} {\mathfrak u}(n) X^n$ 
equals $\prod_{n \geq 0} (1- X^{2^n})$).
$$
(*) \ \ \ X^a S_a(X) = \frac{-1}{1-X}\left(\sum_{0 \leq m \leq a-1} {\mathfrak u}(m) X^m\right) 
+ (1+X^a) \sum_{n \geq 0} {\mathfrak u}(n) X^{2n}.
$$
Now recall that $a$ is odd, say $a = 2b+1$. Then,
$$
\begin{array}{lll}
\displaystyle\sum_{0 \leq m \leq a-1} {\mathfrak u}(m) X^m &=& 
\displaystyle\sum_{0 \leq m \leq 2b} {\mathfrak u}(m) X^m \\
&=&   \displaystyle\sum_{0 \leq m \leq b} {\mathfrak u}(2m) X^{2m}
    + \displaystyle\sum_{0 \leq m \leq b-1} {\mathfrak u}(2m+1) X^{2m+1} \\
&=&   \displaystyle\sum_{0 \leq m \leq b} {\mathfrak u}(m) X^{2m}
    - \displaystyle\sum_{0 \leq m \leq b-1} {\mathfrak u}(m) X^{2m+1} \\
&=& \displaystyle{\mathfrak u}(b) X^{2b} 
    + (1-X)\left(\displaystyle\sum_{0 \leq m \leq b-1} {\mathfrak u}(m) X^{2m}\right). \\
\end{array}
$$
Using $(*)$ we have
$$
X^a S_a(X) = \frac{-1}{1-X}\left(\sum_{0 \leq m \leq 2b} {\mathfrak u}(m) X^m\right) 
+ (1+X^{2b+1}) \sum_{n \geq 0} {\mathfrak u}(n) X^{2n}
$$
hence
$$ 
X^a S_a(X) = - \frac{{\mathfrak u}(b) X^{2b}}{1-X} 
             - \sum_{0 \leq m \leq b-1} {\mathfrak u}(m) X^{2m}
             + (1+X^{2b+1}) \sum_{n \geq 0} {\mathfrak u}(n) X^{2n}. 
$$
Thus
\begin{multline*}
(**) \ \ X^a S_a(X) = - {\mathfrak u}(b) X^{2b} \sum_{j \geq 0} X^j
                      - \displaystyle\sum_{0 \leq m \leq b-1} {\mathfrak u}(m) X^{2m} \\
                               + (1+X^{2b+1})\displaystyle\sum_{n \geq 0} {\mathfrak u}(n) X^{2n}.
\end{multline*}
In order to finish the proof of the case $a$ odd, we have to prove that all the coefficients of 
$S_a(X)$ either take only the values $0$, $2$, or take only the values $0$, $-2$. Of course it 
suffices to prove the same claim for the coefficients of $X^a S_a(X)$. Write $X^a S_a(X) = \sum c_k X^k$.
Looking at $(**)$ we see that
\begin{itemize}
\item for $k \leq a-1$, we have clearly $c_k = 0$;
\item for $2k \geq a$, we have 
          $c_{2k} = - {\mathfrak u}(b) + {\mathfrak u}(k) \in \{- {\mathfrak u}(b) \pm 1\}$;
\item for $2k+1 \geq a$, we have 
          $c_{2k+1} = - {\mathfrak u}(b) + {\mathfrak u}(k-b) \in \{- {\mathfrak u}(b) \pm 1\}$.
\end{itemize}
Hence, either $u(b) = +1$ and all coefficients of $X^a S_a(X)$ belong to $\{-2, 0 \}$, or
$u(b) = -1$ and all coefficients of $X^a S_a(X)$ belong to $\{0, 2 \}$, and we are done.

\medskip

-- To address the case $a$ where even, it suffices to prove that if the statement in Theorem~\ref{shev} 
is true for some integer $a \geq 1$, then it is true for $2a$. Recall that the Thue-Morse sequence can 
be generated by iteratively applying to $+1$ the morphism $\sigma$ defined on $\{+, -\} := \{+1, -1\}$ 
by $\sigma(+) = + \ -$, $\sigma(-) = - \ +$. Thus the prefix of length $2^{d+1}$ of the Thue-Morse sequence 
is equal to the image by $\sigma$ of the prefix of length $2^d$. 
Hence it suffices to prove that for any $a \geq 1$ one has $C_{2a} = 2C_a \cup (2C_a+1)$. The property 
${\mathfrak u}(2n) = {\mathfrak u}(n)$ and ${\mathfrak u}(2n+1) = - {\mathfrak u}(n)$ then gives us the 
desired conclusion. But $n$ belongs to $C_{2a}$ if and only ${\mathfrak u}(n+2a) = {\mathfrak u}(n)$. 
This happens if and only if either $n=2k$ and ${\mathfrak u}(2k+2a) = {\mathfrak u}(2k)$, or 
$n=2k+1$ and ${\mathfrak u}(2k+1+2a) = {\mathfrak u}(2k+1)$.
This is equivalent to either $n=2k$ and ${\mathfrak u}(k+a) = {\mathfrak u}(k)$, or $n=2k+1$ and 
$-{\mathfrak u}(k+a) = -{\mathfrak u}(k)$. This is exactly saying that $n$ belongs to $C_{2a}$ 
if and only if either $n$ belongs to $2C_a$ or $n$ belongs to $2C_a+1$ (note that these sets are, 
of course, disjoint).

\medskip

-- To finish the proof, we follow exactly the same steps, with ${\mathfrak u}(n+a) + {\mathfrak u}(n)$
replaced by ${\mathfrak u}(n+a) - {\mathfrak u}(n)$. Note that this gives in passing the fact that
$\beta_a = \gamma_a$.
\begin{remark}
Shevelev proved in \cite{Shevelev} his conjecture for the case where $a = 2^r$. The case $a=1$ 
was actually proven by Bernhardt \cite{Bernhardt} (see also \cite{BMRSW}). It was also given
by P. Del\'eham (see his comment dated March 16 2004 in \cite[sequence A003159]{oeis}). Note that this 
case is related to the properties of the period-doubling sequence $({\mathfrak z}(n))_{n \geq 0}$.
The period-doubling sequence is defined as the infinite fixed point of the morphism defined on
$\{+, - \} = \{+1, -1\}$ by $- \to - \ +$, $+ \to - \ -$. This morphism occurs, in particular, in the 
study of iterations of unimodal continuous functions. It is not difficult to see that 
${\mathfrak u}(n) = \prod_{0 \leq k \leq n-1} {\mathfrak z}(k)$ (where, as usual, an empty 
product is equal to $+1$). Thus $C_1$ is the set of $n$'s such that ${\mathfrak z}(n) = 1$.

\medskip

The first values of the sequences $B_a$ and $C_a$ for small values of $a$'s are given 
in the On-Line Encyclopedia of Integer Sequences \cite{oeis}, e.g.,
$C_1=$ A079523; $C_2=$ A081706; $C_3=$ A161579; $C_4=$ A161627; $C_5=$ A161817;
$C_6=$ A161824; $C_7=$ A162311; $C_8=$ A161639; $C_9=$ A161890. 

\medskip

There is a straightforward generalization of Shevelev's question. Let $(z(n))_{n \geq 0}$ be
a $\pm 1$ sequence taking each of the values $\pm 1$ infinitely often. For each integer 
$a \geq 1$, let $B_a := \{\ell_0 < \ell_1 < \ell_2 < \ldots\}$ 
(resp.\ $C_a = \{m_0 < m_1 < m_2 < \ldots\}$) denote the sequence of integers in increasing 
order satisfying $z(\ell + a) = - z(\ell)$ (resp.\ $z(m+a) = z(m)$). These two sequences 
of integers reflect how different sequence $u$ and each of its shifted sequences are. Now 
evaluate this difference in terms of the sequence $u$ iteself: let $(\beta_a(n))_{n \geq 0}$ 
(resp.\ $(\gamma_a(n))_{n \geq 0}$) be the $\pm 1$ sequences defined by $\beta_a(n) := z(\ell_n)$ 
(resp.\ $\gamma_a(n) := z(m_n)$. What can be said about sequences $\beta_a$ and $\gamma_a$?
Or for which sequences $u$ do the sequences $\beta_a$ and $\gamma_a$ have interesting properties?
\end{remark}

\subsection{The zeros of the Thue-Morse Dirichlet series}

In a 1985 paper \cite{AllCohen} H. Cohen and the author looked at the Dirichlet series
$\sum_{n \geq 0} \frac{{\mathfrak u}(n)}{(n+1)^s}$ a priori defined for $\Re s > 1$. They proved 
in particular that this series admits an analytic continuation to the whole plane. 
This continuation admits ``trivial zeros'' (all non-positive integers), and 
``non-trivial zeros'' (the complex numbers $2ik\pi/\log 2$ for $k$ integer). 
This suggests the following Riemann-like hypothesis.

\begin{question}
Is it true that $\sum_{n \geq 0} \frac{{\mathfrak u}(n)}{(n+1)^s}$ has no other non-trivial zeros 
than the complex numbers $2ik\pi/\log 2$ for $k$ integer?
\end{question}

\begin{remark}
It might well be that this question is very difficult to answer. It might also well be the case
that the answer would have no interest at all: a proof of the Riemann hypothesis for the zeta 
function would give precious information on the prime numbers, but we do not know of any 
consequence of the Riemann-like hypothesis above.
\end{remark}

\subsection{Looking for a ``simple'' expression for a certain infinite product}

The following infinite product and its simple expression were given in \cite{Woods, Robbins}
$$
P := \prod_{n \geq 0} \left(\frac{2n+1}{2n+2}\right)^{{\mathfrak u}(n)} = \frac{\sqrt{2}}{2}\cdot
$$
While several generalizations can be found in the literature 
(see, e.g., \cite{AllCohen, ASlms, AllSon, LouPro}), no ``simple'' value
is known for a very similar product $Q$ given below. Namely a strange product 
appears in \cite[p.~193]{FlaMar}:
$$
R := \prod_{n \geq 1} \left(\frac{(4n+1)(4n+2)}{4n(4n+3)}\right)^{{\mathfrak u}(n)}.
$$
Actually it is easily proven that $R=\frac{3}{2Q}$, where
$$
Q := \prod_{n \geq 1} \left(\frac{2n}{2n+1}\right)^{{\mathfrak u}(n)}.
$$

\begin{question}
Does the infinite product $Q$ above have a ``simple'' value? Is it a transcendental number?
Is the Flajolet-Martin constant $\varphi$ defined by 
$\varphi := 2^{-1/2}e^{\gamma}\frac{2}{3} R = 0.77351...$ transcendental?
\end{question}

\begin{remark}
The author gave an easy proof for the value of $P$ (which was written down in \cite{AS} 
and \cite{AllSon}) by computing the product $PQ$. We do not resist giving an even 
easier proof. Using the relations ${\mathfrak u}(2n) = {\mathfrak u}(n)$ and 
${\mathfrak u}(2n+1) = - {\mathfrak u}(n)$, we can write 
$$
\begin{array}{lllll}
P &=& \displaystyle\prod_{n \geq 0}\left(\frac{2n+1}{2n+2}\right)^{{\mathfrak u}(n)} 
&=& \displaystyle\left(\prod_{n \geq 0}\left(\frac{2n+1}{2n+2}\right)^{{\mathfrak u}(2n+1)}\right)^{-1} \\
&=& \displaystyle\left(\frac{\displaystyle\prod_{n \geq 1}\left(\frac{n}{n+1}\right)^{{\mathfrak u}(n)}}
                        {\displaystyle\prod_{n \geq 1}\left(\frac{2n}{2n+1}\right)^{{\mathfrak u}(2n)}}    
             \right)^{-1} 
&=& \left(\displaystyle\frac{\displaystyle\prod_{n \geq 1}\left(\frac{n}{n+1}\right)^{{\mathfrak u}(n)}}
                        {\displaystyle\prod_{n \geq 1}\left(\frac{2n}{2n+1}\right)^{{\mathfrak u}(n)}}    
             \right)^{-1} \\
&=& \left(
          \displaystyle\prod_{n \geq 1}\left(\frac{2n+1}{2(n+1)}\right)^{{\mathfrak u}(n)}\right)^{-1}. \\
\end{array}
$$
Hence, multiplying by $P$,
$$
P^2 = \left(\prod_{n \geq 1}\left(\frac{2n+1}{2(n+1)}\right)^{{\mathfrak u}(n)}\right)^{-1}
\prod_{n \geq 0}\left(\frac{2n+1}{2n+2}\right)^{{\mathfrak u}(n)} = \frac{1}{2}\cdot
$$ 
Thus
$$
\displaystyle\prod_{n \geq 0} \left(\frac{2n+1}{2n+2}\right)^{{\mathfrak u}(n)} = \frac{1}{\sqrt{2}}\cdot
$$
\end{remark}

\subsection{Looking for ``another'' proof of Cobham's theorem}

Cobham's theorem asserts that a sequence which is both $q$- and $r$-automatic, where
$q$ and $r$ are multiplicatively independent (i.e., $\log q / \log r$ irrational), must 
be ultimately periodic \cite{Cobham}. The proof of Cobham is very technical, as are
more recent proofs. The quest for a ``simple'' --or at least more ``conceptual''--
proof whose early start was \cite{Hansel} is still open.

\subsection{Almost everywhere automatic sequences}

Deshouillers, whose motivation was the results in \cite{Des1, Des2}, asked the following 
question \cite{Des3}.  Call a sequence $(v(n))_{n \geq 0}$ {\it almost everywhere 
$d$-automatic} if there exists a $d$-automatic sequence $(w(n))_{n \geq 0}$ such that 
$v$ and $w$ are equal almost everywhere (i.e., such that the set of $n$ for which 
$u(n) \neq v(n)$ has natural density $0$, i.e., 
$\displaystyle\lim_{N \to \infty} \frac{1}{N}\sharp\{n \leq N, \ v_n \neq w_n\} = 0$).
Define similarly {\it almost everywhere periodic} sequences, {\it almost everywhere ultimately
periodic sequences}, and {\it almost everywhere constant} sequences. A simple example is the 
characteristic function of squares: this sequence is almost everywhere constant (actually almost
everywhere equal to $0$).

\bigskip

\begin{question} Is the following generalization of Cobham's theorem true? Let $d_1$ and $d_2$ 
be two multiplicatively independent integers (i.e., $\log d_1/\log d_2$ is irrational). If the 
sequence $(v(n))_{n \geq 0}$ is both almost everywhere $d_1$-automatic and almost everywhere 
$d_2$-automatic, then it is almost everywhere ultimately periodic.
\end{question}

\subsection{Sum of digits, pseudo-randomness, distribution modulo $1$}

The Thue-Morse sequence or its version not reduced modulo $2$ (i.e., the sum of binary digits
of the integers) was used in several questions related to pseudo-randomness (in the sense of Bass
or in the sense of Bertrandias), see, e.g., the paper of Mend\`es France (see \cite{Mendes} and
the references therein), or to distribution modulo $1$. The subject continues to be explored, see, 
e.g., \cite{Abou-Liardet}. A nice question of Gelfond was recently answered in \cite{Mauduit-Rivat}.

\begin{theorem}[Mauduit-Rivat]
The sum of $q$-ary digits of the prime numbers is uniformly distributed 
in the non-trivial arithmetic progressions.
\end{theorem}

\subsection{Algebraic independence of power series on a finite field}

Let $p$ be a prime number. Christol's theorem (see \cite{Christol, CKMFR}) asserts that the 
formal power series $\sum a_n X^n$ in ${\mathbb Z}/p{\mathbb Z}[[X]]$ is algebraic over
${\mathbb Z}/p{\mathbb Z}(X)$ if and only if the $p$-kernel of the sequence $(a_n)_{n \geq 0}$,
i.e., the set of subsequences $\{(a_{p^k n + r)_{n \geq 0}}, \ k \geq 0, \ r \in [0, p^k - 1]\}$ 
is finite. In other words a combinatorial property of the sequence $(a_n)_{n \geq 0}$ is
equivalent to the algebraicity of the associated formal power series. A tempting question is
then whether there exists some combinatorial property of two sequences $(a_n)_{n \geq 0}$ and
$(b_n)_{n \geq 0}$ that is equivalent to the property that the formal power series $\sum a_n X^n$ 
and $\sum b_n X^n$ are algebraically dependent over ${\mathbb Z}/p{\mathbb Z}(X)$. Such a condition
could be the finiteness of some ``kernel'' of sequence $(a_n)_{n \geq 0}$, where the extracted
subsequences would somehow depend on sequence $(b_n)_{n \geq 0}$. G. Christol told us (private
communication) that he does not believe that such a condition can exist in general. If so, we might
be able to find subcases where it would be possible to find such a combinatorial condition?

\subsection{$D$-finite formal power series and automatic sequences}

The formal power series $\sum_{n \geq 0} {\mathfrak t}(n) X^n$ considered as an element of
${\mathbb Q}[[X]]$ is transcendental over ${\mathbb Q}(X)$. This is, e.g., a consequence of
a theorem of Fatou \cite{Fatou}: {\it a power series $\sum_{n \geq 0} a_n z^n$ 
with integer coefficients that converges inside the unit disk is either rational or 
transcendental over ${\mathbb Q}(X)$.} Actually more is known: the formal power series 
$\sum_{n \geq 0} {\mathfrak t}(n) X^n$ being irrational and having its radius of convergence
equal to $1$, it admits the unit circle as natural boundary (from a theorem of Carlson \cite{Carlson} 
extending Fatou's).
Another proof of this transcendence using Christol's theorem can be given. If the formal 
power were algebraic, then the series $\sum_{n \geq 0} ({\mathbf t}_n \bmod 2) X^n$
(resp.\ the series $\sum_{n \geq 0} ({\mathbf t}_n \bmod 3)X^n$) would be algebraic over
${\mathbb F}_2(X)$ (resp.\ ${\mathbb F}_3(X)$). Hence, the sequence $({\mathbf t}_n)_{n \geq 0}$
with values $0, 1$ would be both $2$-automatic and $3$-automatic. From a theorem of Cobham 
\cite{Cobham} the sequence $({\mathbf t}_n)_{n \geq 0}$ would thus be ultimately periodic which is 
not true (recall that $({\mathbf t}_n)_{n \geq 0}$ does not contain cubes).

\bigskip

More generally, if a formal power series $\sum a_n X^n$ has all its coefficients in the field of 
rational numbers ${\mathbb Q}$, and is algebraic of degree $d$ over ${\mathbb Q}(X)$, we know, as 
a consequence of Eisenstein's theorem (announced by Eisenstein in \cite{Eis}, and proved by Heine 
in \cite{Hei}) that the set of prime numbers ${\mathcal P}$ that divide the denominator of at least 
one coefficient $a_n$ is finite. Hence it makes sense to reduce the series modulo any prime not in 
${\mathcal P}$. The reduction modulo such a prime $p$ is clearly algebraic over the field 
${\mathbb F}_p(X)$ of degree $d_p \leq d$. Hence from Christol's theorem, the sequence 
$(a_n \bmod p)_{n \geq 0}$ is $p$-automatic. A strategy for proving that the series $\sum a_n X^n$ 
is transcendental over ${\mathbb Q}(X)$ can thus be to prove that either it is algebraic of 
degree $d_p$ for all but finitely many primes $p$, but that the $d_p$'s are unbounded, or to prove that 
for some prime $p$ that does not divide any of the denominators of the $a_n$'s the sequence 
$(a_n \bmod p)_{n \geq 0}$ is not $p$-automatic. Several examples of this strategy are given, 
e.g., in \cite{All-trans}.

\bigskip

Now, a notion generalizing algebraicity for power series is the notion of {\it D-finiteness}, also 
called {\it holonomy}. A very good reference is \cite{Stanley}. A formal power series is called
{\it differentiably finite} (D-finite for short) or {\it holonomic} if it satisfies a linear differential 
equations with polynomial coefficients. Examples are the exponential series $f(X) = e^X$ that satisfies
$f' = f$ or the series $g(X) = \log(1+X)$ that satisfies $(1+X)g'(X) = 1$. It is not hard to see that 
any algebraic formal power series is D-finite. This notion makes sense in zero characteristic (in
characteristic $p$ the $p$th derivative of any formal power series is equal to $0$).
A natural question is thus whether there is a way to prove that a formal power series with say integer 
coefficients is not D-finite by proving that its projections modulo prime numbers do or do not 
satisfy some specific properties? The following question is a conjecture due to Christol \cite{Christol-jap}.

\begin{question} 
It is true that a {\it globally bounded} D-finite formal power series is {\it globally automatic}?
In other words, is it true that a formal power series $\sum a_n X^n$ belonging to ${\mathbb Q}[[X]]$
having a non-zero radius of convergence as a series on ${\mathbb C}$ and for which there exist
$\alpha, \beta$ in ${\mathbb Q}$ such that $\alpha \sum a_n (\beta X)^n$ belongs to ${\mathbb Z}[[X]]$
has the property that for all but finitely many primes $p$ and for all positive integers $h$ the sequence
$(a_n \bmod p^h)_{n \geq 0}$ is is $p$-automatic (or equivalently $p^h$-automatic).
\end{question}

\begin{remark}
As noted by Christol the answer to the above question is yes if two classical conjectures respectively 
due to Bombieri and Dwork are true (see precise formulations in \cite{Christol-jap}).
\end{remark}

\subsection{Transcendence of morphic real numbers and morphic continued fractions}

Several questions can be asked about the ``transcendence'' of the Thue-Morse sequence ${\mathfrak t}$.
For example: 
\begin{itemize}

\item is the real number $\sum_{n \geq 0} {\mathfrak t}(n)/2^n$ transcendental? 
\item is the formal power series $\sum_{n \geq 0} {\mathfrak t}(n) X^n \in {\mathbb Q}[[X]]$ 
      transcendental over ${\mathbb Q}(X)$? 
\item is the formal power series $\sum_{n \geq 0} ({\mathfrak t}(n) \bmod 2) X^n \in {\mathbb Z}/2 {\mathbb Z}[[X]]$
      transcendental over ${\mathbb Z}/2 {\mathbb Z}(X)$? 
\item is the formal power series $\sum_{n \geq 0} ({\mathfrak t}(n) \bmod 3) X^n \in {\mathbb Z}/3 {\mathbb Z}[[X]]$
      transcendental over ${\mathbb Z}/3 {\mathbb Z}(X)$? 
\item is the continued fraction $[1+{\mathfrak t}(0), 1+{\mathfrak t}(1), 1+{\mathfrak t}(2), \ldots]$ 
      transcendental over ${\mathbb Q}$?

\end{itemize}

\noindent
Similar questions can be asked by replacing the Thue-Morse sequence (which, as recalled above, is the 
iterative fixed point beginning with $0$ of the --uniform-- morphism $0 \to 01$, $1 \to 10$) by iterative
fixed points of --not necessarily uniform-- morphisms, or even by morphic sequences, i.e., pointwise images 
of iterative fixed points of morphisms. A survey about these questions, and in particular the fact that
the answers to the questions above about transcendence are respectively ``yes, yes, no, yes, yes'' can
be found in \cite{All-alg}. More results that were either not published, or not proved at the time of that 
survey, are due in particular to Adamczewski and Bugeaud, and can be found in the nice book \cite{Bugeaud}. 
But, though we now know that ``automatic real numbers'' (real numbers whose expansion 
in some integer base is automatic, i.e., is the pointwise image of the iterative fixed point of a uniform 
morphism) are either rational or transcendental, and that ``automatic continued fractions'' (i.e., continued 
fractions whose sequences of partial quotients are automatic) are either quadratic or transcendental,
the general question about transcendence of real numbers whose expansion in some integer base, or whose
continued fraction expansion, is the pointwise image of the iterative fixed point of some general morphism
is still open.

\section{Acknowledgments} We would like to thank V. Shevelev heartily for several interesting email 
discussions about the Thue-Morse sequence and about his conjecture.

\end{document}